\input amstex
\define\sgn{\operatorname{sgn}}
\documentstyle{amsppt}
\pagewidth{12.79 cm}
\topmatter
\title 
Remarks on square functions in the Littlewood-Paley theory
\endtitle
\author
Shuichi Sato
\endauthor
\address
\ 
\newline 
Department of Mathematics 
\newline
Faculty of Education
\newline
Kanazawa University   
\newline
Kanazawa, 920-11
\newline
Japan
\endaddress
\email
shuichi\@kenroku.ipc.kanazawa-u.ac.jp
\endemail
\abstract\nofrills  
We prove that certain square function operators in the Littlewood-Paley theory defined by the kernels 
without  any regularity are bounded on $L^p_w$, $1<p<\infty$, $w\in A_p$ 
 (the weights of Muckenhoupt).   Then, we  give some applications to the Carleson measures on  
 the upper half space. 
\endabstract
\endtopmatter
\document

\head 1.  Introduction \endhead

In this note we shall prove the weighted $L^p$-estimates for the Littlewood-Paley type square 
functions arising from kernels satisfying only size and cancellation conditions.  
Suppose that $\psi \in L^1(\bold R^n)$ satisfies  
$$\int_{\bold R^n} \psi (x)\,dx = 0. \tag 1.1$$  
  We consider a square function of Littlewood-Paley type    
$$S(f)(x) = S_{\psi}(f)(x) = \left( \int_0^{\infty}|\psi_t \star f(x)|^2\,\frac{dt}{t} 
\right)^{1/2},$$ 
where $\psi_t(x) = t^{-n}\psi(t^{-1}x)$. 

If $\psi$  satisfies, in addition to (1.1),  
 
$$|\psi(x)|\leq c(1+|x|)^{-n-\epsilon} \qquad\text{for some \quad $\epsilon >0$} \tag 1.2 $$
$$\int_{\bold R^n}|\psi(x-y) - \psi(x)|\,dx \leq c|y|^\epsilon 
\qquad\text{for some \quad $\epsilon >0$}, \tag 1.3 $$ 
then it is known that the operator $S$ is bounded on $L^p(\bold R^n)$ for all $p\in (1, \infty)$ 
(see Benedek, Calder\'on and Panzone \cite{1}).  Well-known examples are as follows.   
\example{Example 1} Let $P_t(x)$ be the Poisson kernel for the upper half space $\bold R^n \times 
(0, \infty)$: 
$$ P_t(x)=c_n\frac{t}{(|x|^2+t^2)^{(n+1)/2}} .$$
 Put 
$$\psi (x)= \left( \frac{\partial}{\partial t}P_t(x)\right)_{t=1}.$$
Then, $S_{\psi}(f)$ is the Littlewood-Paley $g$ function. 
\endexample
\example{Example 2}
Consider the Haar function $\psi$ on $\bold R$ :  
$$\psi(x) = \chi_{[-1,0]}(x) - \chi_{[0,1]}(x),$$
where $\chi_E$ denotes the characteristic function of a set $E$.
Then, $S_{\psi}(f)$ is the Marcinkiewicz integral 
$$ \mu(f)(x)= \left( \int_0^{\infty}|F(x+t)+F(x-t)-2F(x)|^2\,\frac{dt}{t^3} \right)^{1/2},$$
where $F(x)=\int_0^xf(y)\,dy$.
\endexample

In this note, we shall prove that the $L^p$-boundedness of $S$ still 
holds without the assumption (1.3);  the conditions (1.1) and (1.2) only are sufficient.  
 This has been already known for the $L^2$-case 
(see Coifman and Meyer \cite{3, p. 148}, and also Journ\'e \cite{7, pp. 81-82} for a proof). 

To state our result more precisely,  
we consider the least non-increasing radial majorant of $\psi$
$$h_{\psi}(|x|) = \sup_{|y|\geq |x|}|\psi(y)|.$$
 We also need to consider two seminorms   
$$B_{\epsilon}(\psi) = \int_{|x|>1}|\psi(x)|\,|x|^{\epsilon}\,dx   
\qquad \text{ for \quad  $\epsilon >0$} , $$
$$D_u(\psi)=\left(\int_{|x|<1}|\psi(x)|^u\,dx\right)^{1/u} \qquad \text{for \quad $u>1$}.$$
We shall prove the following result.

\proclaim{Theorem 1} Put $H_\psi(x)=h_{\psi}(|x|)$.  
If $\psi \in L^1(\bold R^n)$ satisfies $(1.1)$ and 
\roster
\item "({1})" $B_{\epsilon}(\psi) < \infty$  \quad  for some $\epsilon >0 $ ; 
\item "({2})"  $D_u(\psi)<\infty$ \quad for some $u>1$ ;
\item "({3})" $H_\psi \in L^1(\bold R^n)$ ;  
\endroster 
 then the operator $S_\psi$ is bounded on $L^p_w$ : 
$$\|S_\psi(f)\|_{L^p_w} \leq C_{p,w}\|f\|_{L^p_w}$$
 for all $p \in (1, \infty)$ and $w \in A_p$, where $A_p$ denotes the weight 
 class of Muckenhoupt (see \cite{6,7}), and 
 $$  \|f\|_{L^p_w}=\|f\|_{L^p(w)}= \left(\int_{\bold R^n} |f(x)|^pw(x)\,dx\right)^{1/p}.$$
\endproclaim  

In fact, we shall prove a more general result.  

\proclaim{Theorem 2} Suppose that $\psi\in L^1(\bold R^n)$ satisfies $(1.1)$ and 
\roster
\item "({1})" $B_{\epsilon}(\psi)<\infty$ \quad for some $\epsilon > 0$ ; 
\item "({2})"  $D_u(\psi)<\infty$ \quad for some $u>1$ ;
\item "({3})" $|\psi(x)| \leq h(|x|)\Omega(x')$ $(x'=|x|^{-1}x)$ for some non-negative functions 
$h$ and $\Omega$  such that   
\item "({a})" $h(r)$ is non-increasing for $r \in (0, \infty)$ ;
\item "({b})"  if $H(x)=h(|x|)$,    $H \in L^1(\bold R^n)$ ;
\item "({c})" $\Omega \in L^q(S^{n-1})$ for some  $q$, $2\leq q \leq \infty$. 
\endroster
Then, the operator $S_\psi$ is bounded on $L^p_w$ for $p>q'$ and $w \in A_{p/q'}$, 
where $q'$ denotes the conjugate exponent of $q$.    
\endproclaim

 When $\psi$ is compactly supported, we have another formulation, which is not included 
 in Theorem 2.    

\proclaim{Theorem 3}  Suppose that $\psi\in L^1(\bold R^n)$ satisfies $(1.1)$ and 
\roster
\item "({1})" $\psi$ is compactly supported ;
\item "({2})" $\psi \in L^q(\bold R^n)$ for some $q\geq 2$.
\endroster 
Then  $S_\psi : L^p_w \to L^p_w$ for $p>q'$ and $w \in A_{p/q'}$.  
\endproclaim

 These results will be derived from more abstract ones.  
Let $\psi \in L^1(\bold R^n)$ satisfy (1.1). 
We also assume the following : 
\roster
\item "({1})" There exists $\epsilon \in (0, 1)$ such that 
$$ \int_1^2|\hat{\psi}(t\xi)|^2\, dt \leq 
c\min\left(|\xi|^{\epsilon}, |\xi|^{-\epsilon}\right) \qquad 
\text{for all \quad $\xi \in \bold R^n$ ,} \tag 1.4 $$
where $\hat{\psi}$ denotes the Fourier transform 
$$\hat{\psi}(\xi) = \int\psi(x)e^{-2\pi i\langle x,\xi\rangle}\,dx ,\qquad 
\langle x,\xi\rangle = \sum_{j=1}^nx_j\xi_j \quad\text{(the inner product in $\bold R^n$)}.$$

\item "({2})" Let $1\leq s \leq 2$.  For all $w \in A_s$,  we have 
$$  \sup_{k \in \bold Z} \int_{\bold R^n} \int_1^2 \left|\psi_{t2^k}\star f(x)\right|^2\,dt\, 
w(x) \,dx 
\leq C_w\|f\|_{L^2_w}^2 \qquad \text{for all \quad $f \in \Cal S(\bold R^n)$ ,}  \tag 1.5 $$
where $\bold Z$ denotes the integer group and $\Cal S(\bold R^n)$ the Schwartz space. 
\endroster 
Under these assumptions  the following holds. 

\proclaim{Proposition 1}  For $p > 2/s$ and $w \in A_{ps/2}$,  the operator $S_\psi$  
  is bounded  on $L^p_w$.  
\endproclaim

This will be used to prove the next result.  

\proclaim{Proposition 2} Put 
$$J_{\epsilon}(\psi)=\sup_{|\xi|=1}\iint_{\bold R^n \times\bold R^n}|\psi(x)\psi(y)|
\left|\langle \xi, x-y\rangle \right|^{-\epsilon}\,dx\,dy 
\qquad\text{for \quad $\epsilon\in (0,1]$}. $$ 
Let $\psi \in L^1$ satisfy $(1.1)$ and $(1.5)$.  Then if 
$B_{\epsilon}(\psi)<\infty$ and  $J_{\epsilon}(\psi)<\infty$ for some $\epsilon \in (0,1]$, 
 the operator $S_\psi$ is bounded on $L^p_w$ for $p > 2/s$ and $w \in A_{ps/2}$. 
\endproclaim

In \S 2, we shall prove Proposition 1  by the method of the proof of 
 Duoandikoetxea and Rubio de Francia \cite{5, Corollary 4.2} and then 
 Proposition 2  by using Proposition 1.   Proposition 2 will be applied  
 to prove Theorems 2 and 3 in \S 3.   Finally, in \S 4, we shall give 
 some applications of Theorem 1 to generalized Marcinkiewicz integrals and 
 the Carleson measures on the upper half space $\bold R^n\times (0, \infty)$.

To conclude this section, we state a result for the $L^2$-case, from which the result of 
Coifman-Meyer mentioned above immediately follows, and an idea of the   
proof will be applied later too (see the proof of Lemma 2).

\proclaim{Proposition 3}Suppose that $\psi \in L^1$ satisfies $(1.1)$.  Let 
$$L(\psi)=\sup_{|\xi|=1}\iint_{\bold R^n \times\bold R^n}|\psi(x)\psi(y)|
\left|\log \left|\langle \xi, x-y\rangle \right|\right| \,dx\,dy .$$
Then, if $L(\psi) <\infty$, the operator $S_{\psi}$ is bounded on $L^2$.        
\endproclaim
\demo{Proof}It is sufficient to show that  
$$\sup_{|\xi|=1}\int_{0}^{\infty} \left|\hat{\psi}(t\xi)\right|^2 \, \frac{dt}{t}<\infty.$$
We write 
$$|\hat{\psi}(t\xi)|^2=\hat{\psi}(t\xi)\overline{\hat{\psi}(t\xi)}
=\iint_{\bold R^n \times \bold R^n}
\psi(x)\overline{\psi(y)}e^{-2\pi it\langle \xi, x-y\rangle}\,dx\,dy ,$$
and so
$$\int_{0}^{\infty} \left|\hat{\psi}(t\xi)\right|^2 \, \frac{dt}{t}
=\lim_{N\to\infty,\epsilon\to 0}\iint\psi(x)\overline{\psi(y)}
\left(\int_\epsilon^Ne^{-2\pi it\langle \xi, x-y\rangle}\frac{dt}{t}\right)dx\,dy.$$

Note that
$$\int_\epsilon^N \left(e^{-2\pi it\langle \xi, x-y\rangle}- \cos(2\pi t)\right)\frac{dt}{t}
\to -\log \left|\langle \xi, x-y\rangle\right| - i\frac{\pi}{2}\sgn \langle \xi, x-y\rangle$$
as $N\to \infty$ and $\epsilon\to 0$, and the integral is bounded, uniformly in $\epsilon$ 
and $N$,  by 
$$c\left(1+\left|\log \left|\langle \xi, x-y\rangle \right|\right|\right).$$
Thus, using (1.1) and the dominated convergence theorem, we get  
$$\int_{0}^{\infty} \left|\hat{\psi}(t\xi)\right|^2 \, \frac{dt}{t}=
\iint \psi(x)\overline{\psi(y)}\left(-\log \left|\langle \xi, x-y\rangle\right| - 
i\frac{\pi}{2}\sgn \langle \xi, x-y\rangle\right)\,dx\,dy.$$ 
This immediately implies the conclusion. 
\enddemo

\remark{Remark }
In the one-dimensional case, it is easy to see that if 
$$\int |\psi(x)|\log(2+|x|)\, dx <\infty \qquad \text{and} \qquad
\int |\psi(x)|\log(2+|\psi(x)|)\, dx <\infty, $$
then $L(\psi)<\infty$, and so $S_\psi : L^2 \to L^2$.   
\endremark

\head 2.  Proofs of Propositions 1 and 2 \endhead 
We use a Littlewood-Paley decomposition.   
Let  $f\in \Cal S(\bold R^n)$, and define 
$$\widehat{\Delta_j(f)}(\xi) = \Psi(2^j\xi)\hat{f}(\xi) \qquad \text{for \quad $j \in \bold Z$},$$ 
where $\Psi \in C^{\infty}$ is supported in $\{1/2\leq |\xi|\leq 2\}$ and satisfies  
$$\sum_{j\in \bold Z} \Psi(2^j\xi) = 1 \qquad \text{for \quad $\xi \neq 0$.}$$ 
 Decompose 
 $$f\star \psi_t(x) = \sum_{j\in \bold Z}\sum_{k\in \bold Z}\Delta_{j+k}(f\star \psi_t)(x)
 \chi_{[2^k, 2^{k+1})}(t) = \sum_{j\in \bold Z}F_j(x, t), \quad \text{say},$$   
 and define 
 $$T_j(f)(x) = \left( \int_0^{\infty}|F_j(x,t)|^2\,\frac{dt}{t} \right)^{1/2}.$$ 
 Then 
 $$S(f)(x) \leq  \sum_{j\in \bold Z}T_j(f)(x).$$
 
 Put $E_j =\{2^{-1-j}\leq |\xi| \leq 2^{1-j}\}$. Then by the Plancherel theorem and (1.4) we have 
 $$\align
 \|T_j(f)\|_2^2 &= \sum_{k\in \bold Z}\int_{\bold R^n}\int_{2^k}^{2^{k+1}}
 \left|\Delta_{j+k}\left(f\star \psi_t \right)(x)\right|^2 \, \frac{dt}{t}\, dx
 \\
 &\leq \sum_{k\in \bold Z} c \int_{E_{j+k}}\left(\int_{2^k}^{2^{k+1}}
 \left|\hat{\psi}(t\xi)\right|^2 \, \frac{dt}{t}\right)\,\left|\hat{f}(\xi)\right|^2 \,d\xi
 \\
 &\leq  \sum_{k\in \bold Z} c\int_{E_{j+k}}\min\left(|2^k\xi|^{\epsilon}, |2^k\xi|^{-\epsilon}\right) 
 \left|\hat{f}(\xi)\right|^2 \,d\xi
 \\
 &\leq c 2^{-\epsilon |j|} \sum_{k\in \bold Z}\int_{E_{j+k}}\left|\hat{f}(\xi)\right|^2 \,d\xi
 \\
 &\leq c 2^{-\epsilon |j|}\|f\|_2^2,
 \endalign
 $$
 where the last inequality holds since the sets $E_j$ are finitely overlapping.  (We denote by  
 $\|\cdot\|_p$ the ordinary $L^p$-norm.)

 On the other hand, for $w \in A_s$ by (1.5)  we see that 
 $$\align 
 \|T_j(f)\|_{L^2_w}^2 &= \sum_{k\in \bold Z}\int_{\bold R^n}\int_{2^k}^{2^{k+1}}
 \left|\Delta_{j+k}(f)\star \psi_t (x)\right|^2 \, \frac{dt}{t}\, w(x)\, dx
 \\
 &\leq \sum_{k\in \bold Z} c\int_{\bold R^n}\left|\Delta_{j+k}(f)(x)\right|^2w(x)\, dx
 \\
 &\leq c\|f\|_{L^2_w}^2,
 \endalign
 $$
 where the last inequality follows from a well-known Littlewood-Paley inequality for $L^2_w$ 
 since $A_s \subset A_2$. 
 
 Interpolating with change of measures between the two estimates above, we get 
 $$\|T_j(f)\|_{L^2(w^u)} \leq c2^{-\epsilon(1-u)|j|/2} \|f\|_{L^2(w^u)}$$
 for $u\in (0, 1)$.  If we choose $u$ (close to 1) so that $w^{1/u} \in A_s$, then 
 from this inequality we get 
  $$\|T_j(f)\|_{L^2_w} \leq c2^{-\epsilon(1-u)|j|/2} \|f\|_{L^2_w} , $$
  and so 
  $$\|S(f)\|_{L^2_w} \leq  \sum_{j\in \bold Z}\|T_j(f)\|_{L^2_w} \leq c\|f\|_{L^2_w}.$$ 
  Thus the extrapolation theorem of Rubio de Francia \cite{8} implies the conclusion.

To derive Proposition 2 from Proposition 1 we prepare the following lemmas.
 
\proclaim{Lemma 1}
If $\psi\in L^1(\bold R^n)$ satisfies $(1.1)$ and $B_{\epsilon}(\psi)<\infty$ for  
$\epsilon \in (0, 1]$,  then 
$$|\hat{\psi}(\xi)| \leq c|\xi|^{\epsilon} \qquad \text{for all \quad $\xi \in \bold R^n$}.$$  
\endproclaim
\demo{Proof} Since $a \leq a^{\epsilon}$ for $a, \epsilon \in (0, 1]$, we see that
$$\align
 |\hat{\psi}(\xi)|= \left|\int \psi(x)\left(e^{-2\pi i\langle x,\xi\rangle}-1\right)\, dx\right| 
&\leq c\int |\psi(x)|\min(1, |\langle x,\xi\rangle|)\, dx 
\\
&\leq c|\xi|^{\epsilon}\int |\psi(x)||x|^{\epsilon}\,dx.
\endalign $$
This completes the proof.
\enddemo

\proclaim{Lemma 2} If $\psi \in L^1(\bold R^n)$ and $J_{\epsilon}(\psi) < \infty$ for  
$\epsilon\in (0,1]$, then 
$$\int_1^2|\hat{\psi}(t\xi)|^2\,dt \leq c|\xi|^{-\epsilon} \qquad 
\text{for all \quad $\xi \in \bold R^n$}.$$ 
\endproclaim
\demo{Proof} As in the proof of Proposition 3, we see that  
$$\int_1^2|\hat{\psi}(t\xi)|^2\,dt=\iint_{\bold R^n \times\bold R^n}\psi(x)\overline{\psi(y)}
\frac{e^{-4\pi i\langle \xi, x-y\rangle}-e^{-2\pi i\langle \xi, x-y\rangle}}
{-2\pi i\langle \xi, x-y\rangle}
\,dx\,dy .$$
Thus 
$$\align 
\int_1^2|\hat{\psi}(t\xi)|^2\,dt &\leq  
c \iint_{\bold R^n \times\bold R^n}\left|\psi(x)\psi(y)\right|\min
\left(1, |\langle \xi, x-y\rangle|^{-1}\right)\,dx\,dy 
\\
&\leq cJ_{\epsilon}(\psi)\,|\xi|^{-\epsilon}.
\endalign $$
This completes the proof. 
\enddemo

Now, we can see that Proposition 1 implies Proposition 2, since the condition (1.4) 
follows from Lemmas 1 and  2.

\head 3.  Proofs of Theorems 2 and 3 \endhead 
To get Theorem 2 from Proposition 2 we need Lemmas 3 and 4 below. 
First, we give a sufficient condition for $J_{\epsilon}(\psi) < \infty$. 
 
\proclaim{Lemma 3} Let $h(r)$, $h\geq 0$,  be a non-increasing function for $r>0$ satisfying  
$H\in L^1(\bold R^n)\cap L^\infty(\bold R^n)$, where $H(x)=h(|x|)$,   
and let $\Omega \in L^v(S^{n-1})$, $v>1$, $\Omega\geq 0$. Suppose that $F$ is a non-negative 
function such that 
$$F(x)\leq h(|x|)\Omega(x') \qquad\text{for \quad $|x|>1$}$$
and  $D_u(F)<\infty $ for  $u>1$.  
Then $J_{\epsilon}(F)<\infty$ if $\epsilon <\min(1/u', 1/v')$.   
\endproclaim
\demo{Proof}
For non-negative functions $f$, $g$ and $\xi\in S^{n-1}$  put 
$$L_\epsilon (f,g;\xi)=\iint_{\bold R^n \times\bold R^n}f(x)g(y)
\left|\langle \xi, x-y\rangle\right|^{-\epsilon} \,dx\,dy.$$
Decompose $F$ as $F=E+G$, where $E(x)=F(x)$ if $|x|<1$ and $E(x)=0$ otherwise. 
Then
$$L_\epsilon (F,F;\xi) = L_\epsilon (E,E;\xi)+2L_\epsilon (E,G;\xi)+L_\epsilon (G,G;\xi).$$
We show that each of $L_\epsilon (E,E;\xi), L_\epsilon (E,G;\xi)$ and  $L_\epsilon (G,G;\xi)$ 
is bounded by a constant independent of $\xi$ if $\epsilon <\min(1/u', 1/v')$.  
  
First, by H\"older's inequality  and a change of variables 
$$L_\epsilon (E,E;\xi) \leq \|E\|_u^{2}\left(\iint_{|x|<1, |y|<1}
\left|x_1-y_1\right|^{-\epsilon u'} \,dx\,dy\right)^{1/u'},$$ 
where we note that $\|E\|_u=D_u(F)$.  

Next, by H\"older's inequality again 
$$L_\epsilon (E,G;\xi) \leq \|E\|_u\left(\int_{|x|<1}
\left(\int_{\bold R^n}G(y)\left|x_1-\langle \xi,y\rangle\right|^{-\epsilon}\,dy\right)^{u'}
\,dx\right)^{1/u'}.$$
For $s>0$, let  
$$I_{\epsilon}(s)=\int_{S^{n-1}}
\left|x_1-\langle \xi, s\omega\rangle\right|^{-\epsilon}\Omega(\omega) \,d\sigma(\omega)$$
for fixed $x_1$ and $\xi$, where $d\sigma$ denotes the Lebesgue surface measure of $S^{n-1}$  
(when $n=1$, let $\sigma(\{1\})=\sigma(\{-1\})=1$).   Then by H\"older's inequality 
$$I_{\epsilon}(s)\leq \left(N_{\epsilon v'}(s)\right)^{1/v'}\|\Omega\|_v ,$$
where
$$N_{\epsilon}(s)=\int_{S^{n-1}}\left|x_1 -s\omega_1\right|^{-\epsilon}\,d\sigma(\omega).$$
Thus, using H\"older's inequality, 
$$\align 
\int_{\bold R^n}G(y)\left|x_1-\langle \xi,y\rangle\right|^{-\epsilon}\,dy
&\leq \int_0^\infty h(s)s^{n-1}I_{\epsilon}(s) \,ds
\\
&\leq \|\Omega\|_v\int_0^\infty h(s)s^{n-1}\left(N_{\epsilon v'}(s)\right)^{1/v'} \,ds
\\
&\leq c\|H\|_1^{1/v}\|\Omega\|_v\left(\int_0^\infty h(s)s^{n-1}N_{\epsilon v'}(s) \,ds\right)^{1/v'}
\\
&= c\|H\|_1^{1/v}\|\Omega\|_v
\left(\int_{\bold R^n}h(|y|)\left|x_1-y_1\right|^{-\epsilon v'}\,dy\right)^{1/v'}. 
\endalign$$

Therefore, the desired estimate for $L_\epsilon (E,G;\xi)$ follows if we show that 
$$\sup_{x_1\in \bold R}\int_{\bold R^n}h(|y|)\left|x_1-y_1\right|^{-\epsilon v'}\,dy <\infty .
\tag 4.1 $$
To see this, we split the domain of the integration as follows :
$$ \align   
\int_{\bold R^n}h(|y|)\left|x_1-y_1\right|^{-\epsilon v'}\,dy 
&=\int_{|x_1-y_1|<1}h(|y|)\left|x_1-y_1\right|^{-\epsilon v'}\,dy 
 \\
 &\quad + \int_{|x_1-y_1|>1}h(|y|)\left|x_1-y_1\right|^{-\epsilon v'}\,dy 
 \\
 &=I_1+I_2, \quad \text{say}.
\endalign $$

Clearly $I_2\leq \|H\|_1$.  To estimate $I_1$ we may assume that $n\geq 2$;  
the case $n=1$ can be easily disposed of since $h$ is bounded.   
 We need further splitting of the domain of the integration. 
We write $y=(y_1,y')$, $y'\in \bold R^{n-1}$.   Then 
$$\align   
I_1 &=\int\limits_{\Sb |x_1-y_1|<1\\ |y'|<1\endSb} h(|y|)\left|x_1-y_1\right|^{-\epsilon v'}\,dy 
+ \int\limits_{\Sb |x_1-y_1|<1\\ |y'|>1\endSb} h(|y|)\left|x_1-y_1\right|^{-\epsilon v'}\,dy
\\ 
&=I_3+I_4, \quad\text{say}.
\endalign $$
It is easy to see that   
$$I_3\leq  \|H\|_\infty\,\int_{|y|<2}|y_1|^{-\epsilon v'}\,dy  <\infty .$$
 Next,  since $h(|y|)\leq h(|y'|)$, 
$$\align  
I_4 &\leq \int_{|y_1|<1}|y_1|^{-\epsilon v'}\,dy_1\,\int_{|y'|>1}h(|y'|)\,dy' 
\\
&\leq c\int_{|y_1|<1}|y_1|^{-\epsilon v'}\,dy_1\,\int_{|y|>1}h(|y|)\,dy <\infty .
\endalign $$

It remains  to estimate $L_\epsilon (G,G;\xi)$.   Note that 
$$L_\epsilon (G,G;\xi)\leq \int_0^\infty\!\int_0^\infty h(r)h(s)r^{n-1}s^{n-1}I_{\epsilon}(r, s)
\,dr\,ds,\tag 4.2 $$
where
$$I_{\epsilon}(r,s)=\iint_{S^{n-1}\times S^{n-1}}
\left|\langle \xi, r\theta -s\omega\rangle\right|^{-\epsilon}
\Omega(\theta)\Omega(\omega) \,d\sigma(\theta)\,d\sigma(\omega).$$
 By H\"older's inequality 
$$I_{\epsilon}(r,s)\leq \left(N_{\epsilon v'}(r,s)\right)^{1/v'}\|\Omega\|_v^2, \tag 4.3 $$
where
$$N_{\epsilon}(r,s)=\iint_{S^{n-1}\times S^{n-1}}\left|r\theta_1 -s\omega_1\right|^{-\epsilon}
 \,d\sigma(\theta)\,d\sigma(\omega).$$
Using the estimate (4.3) in (4.2) and then applying H\"older's inequality, we see that 
$$\align 
L_\epsilon (G,G;\xi) &\leq c\|H\|_1^{2/v}\|\Omega\|_v^2
\left(\int_0^\infty\!\!\int_0^\infty N_{\epsilon v'}(r,s)
h(r)h(s)r^{n-1}s^{n-1} \,dr\,ds\right)^{1/v'}
\\
&= c\|H\|_1^{2/v}\|\Omega\|_v^2\left(\iint_{\bold R^n \times\bold R^n}
h(|x|)h(|y|)\left|x_1-y_1\right|^{-\epsilon v'}\,dx\,dy\right)^{1/v'}.
\endalign$$
Therefore, the desired estimates  follows again from (4.1).
This completes the proof.
\enddemo 

For a non-negative function $\Omega$ on $S^{n-1}$ we define 
a non-isotropic Hardy-Littlewood maximal function 
$$M_{\Omega}(f)(x)= \sup_{r>0}r^{-n}\int_{|y|<r}|f(x-y)|\Omega(|y|^{-1}y)\,dy.$$

To prove Theorem 2 we also need the following (see Duoandikoetxea \cite{4}). 

\proclaim{Lemma 4} If $\Omega \in  L^q(S^{n-1})$, $q\geq 2$, and $w\in A_{2/q'}$, then 
$M_{\Omega}$ is bounded on $L^2_w$.   
\endproclaim

Now we can prove Theorem 2. As in Stein \cite{10, pp. 63-64}, we can show that 
$$\sup_{t>0}\left|\psi_{t}\star f(x)\right| \leq c\,M_{\Omega}(f).$$ 
So, by Lemma 4 we see that the condition (1.5) holds for $\psi$ of Theorem 2 
with $s=2/q'$.  

Next, applying Lemma 3, we see that $J_\epsilon(\psi)<\infty$ for $\epsilon < 
\min(1/u', 1/q')$ (note that $h(r)$ of Theorem 2 (3) is bounded for $r\geq 1$).  
Combining these facts with the assumption in Theorem 2 (1), we can apply Proposition 2 
to reach the conclusion. 

Finally, we give the proof of Theorem 3.  Clearly $B_1(\psi)<\infty$, and   
$J_{1/(2q')}(\psi) <\infty$  by applying Lemma 3 suitably. 
Therefore, the conclusion follows from Proposition 2 if we show that the condition (1.5) 
holds with $s=2/q'$. But, for $q>2$ this is a consequence of the inequality 
$$\sup_{t>0}\left|\psi_{t}\star f(x)\right| \leq c\,M(|f|^{q'})^{1/q'},$$ 
where $M$ denotes the Hardy-Littlewood maximal operator.  (This inequality is easily seen 
by H\"older's inequality.) 

To prove the condition (1.5) when $q=2$ and $w\in A_1$, we may assume that $\psi$ is supported 
in $\{|x|<1\}$.  Then by  Schwarz's inequality 
$$\left|\psi_{t}\star f(x)\right|^2\leq t^{-n}\|\psi\|_2^2\int_{|y|<t}|f(x-y)|^2\,dy .$$
Integrating with the measure $w(x)\,dx$ and using a property of the $A_1$-weight function, 
we get 
$$\align 
\int \left|\psi_{t}\star f(x)\right|^2w(x)\,dx &\leq \|\psi\|_2^2
\int |f(y)|^2t^{-n}\int_{|x-y|<t}w(x)\,dx \,dy 
\\
&\leq C_w \|\psi\|_2^2 \int |f(y)|^2w(y)\,dy 
\endalign $$
uniformly in $t$.  From this the desired inequality follows.

\head 4.  Applications \endhead 
It is to be noted that Theorem 1 can be applied to study the 
$L^p_w$-boundedness of generalized Marcinkiewicz integrals.   

\proclaim{Corollary 1} 
For $\epsilon>0$, let 
$$\psi(x)= |x|^{-n + \epsilon}\Omega(x')\chi_{(0, 1]}(|x|)   ,$$
where $\Omega\in L^\infty(S^{n-1})$ and $\int \Omega (x')\,d\sigma(x')=0$. 
Define a Marcinkiewicz integral  
$$\mu(f)(x) = \left(\int_0^\infty\left|\psi_t\star f(x)\right|^2\,\frac{dt}{t}
\right)^{1/2}.$$
Then, the operator $\mu$ is bounded on $L^p_w$ for all $p\in (1, \infty)$ and $w\in A_p$ : 
$$\|\mu(f)\|_{L^p_w}\leq C_{p,w}\|f\|_{L^p_w} .$$
\endproclaim 

This result, in particular, removes the Lipschitz condition assumed for $\Omega$ in 
 Stein \cite{9, {\smc Theorem} 1 (2)}. 

Next, we consider applications to the Carleson measures on the upper half spaces.   
  
\proclaim{Corollary 2} 
Suppose $\psi \in L^1$ satisfies $(1.1)$ and 
$$|\psi(x)| \leq c(1+|x|)^{-n-\epsilon} \qquad \text{for some \quad $\epsilon>0$}.$$
Take $b\in BMO$ and $w \in A_2$.  Then the measure 
$$d\nu(x,t) = \left|\psi_t\star b(x)\right|^2\,\frac{dt}{t}\,w(x)\,dx$$
on the upper half space $\bold R^n\times (0, \infty)$ is a Carleson measure with respect to
 the measure $w(x)\,dx$, that is, 
 $$\nu(S(Q)) \leq C_w\|b\|_{BMO}^2\int_Qw(x)\,dx$$
 for all cubes $Q$ in $\bold R^n$, where 
 $$S(Q)= \{(x, t)\in \bold R^n\times (0, \infty) : x\in Q, 0<t\leq \ell(Q)\},$$
 with $\ell(Q)$ denoting  sidelength of $Q$. 
 \endproclaim
 
 This can be proved by using $L^2_w$-boundedness of the operator  $S_{\psi}$ 
 (see Theorem 1)  as in Journ\'e  
 \cite{7, Chap. 6 III, pp. 85--87}.  In  \cite{7}, a similar result  has been proved with 
 an additional assumption on the gradient of $\psi$. 
 
  Arguing as in \cite{7, Chap. 6 III, p. 87}, by Corollary 2 we can get the following. 
 
 \proclaim{Corollary 3 } Let $\psi$ and $b$ be as in Corollary 2. Suppose $\varphi$ satisfies 
 $$|\varphi(x)| \leq c(1+|x|)^{-n-\delta} $$
 for $\delta>0$. Then, the sublinear operator 
 $$T_b(f)(x)=  \left(\int_0^\infty\left|\psi_t\star b(x)\right|^2
 \left|\varphi_t\star f(x)\right|^2\,\frac{dt}{t}\right)^{1/2}$$
 is bounded on $L^p_w$ for all $p\in (1, \infty)$ and $w\in A_p$ : 
 $$\|T_b(f)\|_{L^p_w}\leq C_{p,w}\|b\|_{BMO}\|f\|_{L^p_w}.$$ 
 \endproclaim

 Here again we don't need the assumption on the gradient of $\psi$.  See Coifman and Meyer 
 \cite{3, p. 149} for the $L^2$-case.

  \proclaim{Corollary 4 } Suppose $\eta \in L^1(\bold R^n)$ satisfies the assumptions of 
  Theorem 1 for $\psi$. 
   Let $\psi$, $\varphi$ and $b$ be as in Corollary 3, and define 
  a paraproduct 
  $$\pi_b(f)(x)= \int_0^\infty \eta_t\star\left((\psi_t\star b)\,(\varphi_t\star f)\right)(x)
  \,\frac{dt}{t}.$$
  Then, the operator $\pi_b$  is bounded on $L^p_w$ for all $p\in (1, \infty)$ and $w\in A_p$ : 
 $$\|\pi_b(f)\|_{L^p_w}\leq C_{p,w}\|b\|_{BMO}\|f\|_{L^p_w}.$$
 \endproclaim
 \demo{Proof} Let $g\in L^2(w^{-1})$, $w\in A_2$. Then, since $w^{-1}\in A_2$, 
  by  Schwarz's inequality, Theorem 1 and Corollary 3, for $0<u<v$, we see that   
 $$\align 
 &\left|\int \int_u^v \eta_t\star\left((\psi_t\star b)\,(\varphi_t\star f)\right)(x)
  \,\frac{dt}{t} g(x)\,dx\right| 
  \\
  &\leq \left(\int \int_u^v\left|\tilde{\eta}_t\star g(x)\right|^2
  \,\frac{dt}{t} \,w^{-1}(x)\, dx \right)^{1/2}\|T_b(f)\|_{L^2(w)}
 \\  
 &\leq C_w\|b\|_{BMO}\|g\|_{L^2(w^{-1})}\|f\|_{L^2(w)}, 
 \endalign$$
 where $\tilde{\eta}(x)=\eta(-x)$. From this estimate we can see that $\pi_b(f)$ is 
  well-defined (see Christ \cite{2, III, \S 3}). 
 Taking the supremum over $g$ with $\|g\|_{L^2(w^{-1})}\leq 1$, we get the $L^2_w$-boundedness, 
 and so  the extrapolation theorem of Rubio de Francia implies the conclusion. 
 This completes the proof.
 \enddemo
 
 See Coifman and Meyer \cite{3, p. 149, PROPOSITION 1} for a similar result in 
 the $L^2$-case.

\enddocument